\documentstyle[12pt,amssymb]{article}
\newtheorem{theorem}{Theorem}[section]
\newtheorem{lemma}{Lemma}[section]

\begin{document}
\title{A CHARACTERIZATION OF SCHMIDT GROUPS BY THEIR
ENDOMORPHISM SEMIGROUPS}
\author{Peeter Puusemp\\
Department of Mathematics\\
Tallinn Technical University\\
Ehitajate tee 5\\
Tallinn 19086, Estonia\\
e-mail: puusemp@edu.ttu.ee} \maketitle

\begin{abstract}
A {\it Schmidt group} is a non-nilpotent finite group in which each proper
subgroup is nilpotent. Each Schmidt group $G$ can be described by three
parameters $p,\, q$ and $v$, where $p$ and $q$ are different primes and $v$
is a natural number, $v\ge 1$. Denote by ${\mathcal{S}}$ the class of all
Schmidt groups which have the same parameters $p,\, q$ and $v$. It is shown
in this paper that the class ${\mathcal{S}}$ can be characterized by the
properties of the endomorphism semigroups of the groups of this class. It
follows from this characterization that if $G\in {\mathcal{S}}$ and $H$ is
another group such that the endomorphism semigroups of $G$ and $H$ are
isomorphic, then $H\in {\mathcal{S}}$, too.
\end{abstract}

\medskip
\noindent {\bf Keywords:} group, Schmidt group, endomorphism, endomorphism
semigroup

\medskip
\noindent {\bf 2000 MSC numbers:} 20F99, 20M20

\subsection*{1. Introduction}

\par The question, whether two algebraic structures are isomorphic when
their automorphism groups or endomorphism semigroups are isomorphic, is a
fundamental but generally difficult problem in algebra. For example,
H.Leptin \cite{hLe} showed in 1960 that if $p\ge 5$, then abelian
$p$-groups are determined by their automorphism groups. Only after 25 years
W.Liebert \cite{wLi} proved that it is so if $p=3$. Finally, in 1998
P.Schultz \cite{pS98} showed that Leptin's result holds for $p=2$, too. We
have found many examples of groups which are determined by their
endomorphism semigroups in the class of all groups. Some of such groups
are: finite abelian groups (\cite{pP75}, Theorem 4.2), non-torsion
divisible abelian groups (\cite{pP99}, Theorem 1), generalized quaternion
groups (\cite{pP76}, Corollary 1). In this paper we give a characterization
of Schmidt groups by their endomorphism semigroups.
\par A {\it Schmidt group} is a non-nilpotent finite group in which each
proper subgroup is nilpotent. This notion is dedicated to O.J.Schmidt, who
first described some properties of the mentioned groups (\cite{oS24}). The
structure of Schmidt groups is well-known (L.Redei \cite{lR47,lR56}). Each
Schmidt group $G$ can be described by three natural numbers $p,\, q$ and
$v$, where $p$ and $q$ are different primes and $v\ge 1$. These numbers are
called {\it parameters of} $G$. A Schmidt group is not uniquely determined
by its parameters. But there exists only one Schmidt group $G$ (up to an
isomorphism) with parameters $p,\, q$ and $v$ such that $G$ is a group of
Miller and Moreno. Let us call it {\it the group of Miller and Moreno with
parameters $p,\, q$ and $v$}. Remember that a group of Miller and Moreno is
a non-abelian finite group in which each proper subgroup is abelian.

Fix parameters $p,\, q$ and $v$ and denote by ${\mathcal{S}}$ the class of
all Schmidt groups which have these parameters. In this paper we shall show
that the class ${\mathcal{S}}$ can be characterized by the properties of the
endomorphism semigroups of the groups of this class (Theorem 3.1). It
follows from this characterization that if $G$ and $H$ are groups such that
$G\in {\mathcal{S}}$ and the endomorphism semigroups $\hbox{End} (G)$ and
$\hbox{End} (H)$ are isomorphic, then $H\in {\mathcal{S}}$ (Theorem 3.2).
Let $u$ be the order of $p$ in the group of units of the residue-class ring
$Z_{q}$ modulo $q$. As corollary from the previous two theorems, it follows
that the group of Miller and Moreno with parameters $p,\, q$ and $v$ is
determined by its endomorphism semigroup in the class of all groups if and
only if $u$ is odd (Theorem 3.3). Remark that the endomorphism semigroups
of the groups of the class ${\mathcal{S}}$ were found in \cite{pP01}.

We shall use the following notations: $G$ - a group; $\hbox{End} (G)$ - the
endomorphism semigroup of the group $G$; $Z(G)$ - center of $G$; $[G:H]$ -
index of the subgroup $H$ in the group $G$; $C_{G}(H)$ - the centralizer of
$H$ in $G$; $N_{G}(H)$ - the normalizer of $H$ in $G$; $I_{0}(G)$ - the set
of all idempotents of $\hbox{End} (G)$ different from $0$ and $1$; $[x]=\{
\, y\in \hbox{End} (G)\ \vert \ y^{2}=y,\ xy=y,\ yx=x\, \} $ ($x\in
\hbox{End} (G)$); $\hat{g} $ - the inner automorphism of $G$, generated by
an element $g\in G$; $p,\, q,\, v$ - parameters of a Schmidt group;
${\mathcal{S}}$ - the class of all Schmidt groups with parameters $p,\, q$
and $v$; $Z_r$ - the residue-class ring modulo $r$; $u$ - the order of $p$
in the group of units of the ring $Z_q$; $Z_r[x]$ - the polynomial ring over
$Z_r$; $[a,\, b]=a^{-1}b^{-1}ab$; $A\, \leftthreetimes \, B$ - the
(internal) semidirect product of $A$ and $B$; $K_{G}(x)=\{ \, y\in
\hbox{End} (G)\ \vert \ yx=xy=y\, \} $; $V_{G}(x)=\{ \, y\in \hbox{Aut}
(G)\ \vert \ yx=x\, \} $; $D_{G}(x)=\{ \, y\in \hbox{Aut} (G)\ \vert \
yx=xy=x\, \} $; $H_{G}(x)=\{ \, y\in \hbox{End} (G)\ \vert \ xy=y,\ yx=0\,
\} $;

Remark that $V_{G}(x),\, D_{G}(x),\, K_{G}(x)$ and $H_{G}(x)$ are
subsemigroups of \break $\hbox{End} (G)$, however, $V_{G}(x)$ and $D_{G}(x)$
are subgroups of $\hbox{Aut} (G)$. Note \break that a map is written in this
paper right from the element on which it acts.

\subsection*{2. Preliminaries}

\setcounter{lemma}{0} \setcounter{section}{2}

Let us cite some easy facts useful in the proofs of the main theorems.

If $x$ is an idempotent of $\hbox{End} (G)$, then $G$ decomposes into the
semidirect product $G=\hbox{Ker} \, x\leftthreetimes \hbox{Im} \, x$.
Clearly, if $x$ and $y$ are two idempotents of $\hbox{End} (G)$, then $x=y$
if and only if $\hbox{Im} \, x=\hbox{Im} \, y$ and $\hbox{Ker} \,
x=\hbox{Ker} \, y$.

\begin{lemma}
\label{L: 2.1} If $x$ is an idempotent of $\hbox{End} (G)$, then the
semigroups $K_{G}(x)$ and $\hbox{End} (\hbox{Im} \, x)$ are isomorphic.
This isomorphism is given by the correspondence $z\longmapsto z\vert
\hbox{Im} \, x$ where $z\in K_{G}(x)$.
\end{lemma}

\begin{lemma}
\label{L: 2.2} If $x\in \hbox{End} (G)$ and $\hbox{Im} \, x$ is abelian,
then $\hat{g} \in V_{G}(x)$ for each $g\in G$.
\end{lemma}

\begin{lemma}
\label{L: 2.3} If $x\in \hbox{End} (G),\ g\in C_{G}(\hbox{Im} \, x)$ and
the subgroup $\hbox{Im} \, x$ of $G$ is abelian, then $\hat{g} \in
D_{G}(x)$.
\end{lemma}

\begin{lemma}
\label{L: 2.4} If $x,\, y\in \hbox{End} (G)$ and $xy=yx$, then $(\hbox{Im}
\, x)y\subset \hbox{Im} \, x$ and $(\hbox{Ker} \, x)y\subset \hbox{Ker} \,
x$.
\end{lemma}

The proofs of lemmas 2.1-2.4 are easy exercises.

In the next part of this section $G$ will everywhere denote an arbitrary
Schmidt group with parameters $p,\, q$ and $v$. By L.Redei
\cite{lR47,lR56},the group $G$ satisfies following properties: 1) $G$
decomposes into a semidirect product $G=G'\leftthreetimes C(q^{v})$ of the
derived subgroup $G'$ of $G$ and a cyclic subgroup $C(q^{v})=<b>$ of the
order $q^{v}$; 2) the derived subgroup $G'$ is a $p$-group; 3) the second
derived subgroup $G''$ of $G$ is a subgroup of the center $Z(G)$ of $G$; 4)
the factor-group $G'/G''$ is an elementary abelian $p$-group of the order
$p^{u}$ ( the integer $u$ is defined later in this section) 5)
$Z(G/G'')=<b^{q}G''>$.

Denote $G^{*}=G/G''$. Next we remember some properties of $\hbox{End} (G)$
and $\hbox{End} (G^{*})$ which are proved in \cite{pP01}.

The map $^{*}:\hbox{End} (G)\longrightarrow \hbox{End} (G^{*}),\ \, \tau
\mapsto \tau ^{*} $, where $(gG'')\tau ^{*}=(g\tau)G''\, (g\in G,\ \tau \in
\hbox{End} (G))$, is a monomorphism. This monomorphism induces an
isomorphism between the semigroups $\hbox{End} (G)\setminus \hbox{Aut} (G)$
and $\hbox{End} (G^{*})\setminus \hbox{Aut} (G^{*})$ of all proper
endomorphisms of $G$ and $G^{*}$.

Assume that $\psi (x)$ is an arbitrary irreducible normalized divisor of
the polynomial
$${x^q-1\over x-1} = x^{q-1}+x^{q-2}+\ldots +x+1\in Z_p[x].$$
Denote by $u$ the degree of $\psi (x)$. Otherwise, $u$ is equal to the
order of the element $p$ in the group of units of the ring $Z_{q}$. Remark
that the lowest natural number $n$ such that $\psi (x)$ is a divisor of
$x^n-1$ in $Z_p[x]$ is equal to $q$ (\cite{lR47}, Proposition 6). Denote by
$\overline{Z} _p[x]$ the factor-ring $Z_p[x]/(\psi (x))$ of $Z_p[x]$ over
the principial ideal $(\psi (x))$, generated by $\psi (x)$. Proper
endomorphisms of $G^{*}$ are given as pairs $[n;\, f(x)]$, where $n\in
Z(q^{v})$ and $f(x)\in \overline{Z} _p[x]$. Two pairs $[n_{1};\, f_{1}(x)]$
and $[n_{2};\, f_{2}(x)]$ are equal if and only if 1) $n_{1}=n_{2} (=n)$; 2)
$f_{1}(x)=f_{2}(x)$ or $q$ is a divisor of $n$. Automorphisms of $G^{*}$
are given as triplets $[n;\, a(x);\, b(x)]$, where $n\in Z(q^{v});\ a(x),\,
b(x)\in \overline{Z} _p[x];\ b(x)\ne 0$ and $q$ is not a divisor of $n$. If
two triplets are distinct, then the corresponding automorphisms of $G^{*}$
are distinct, too. The composition rule of elements of $\hbox{End} (G^{*})$
is given by the following equalities:
$$[n;\, f(x)]\cdot [m;\, g(x)]=[nm;\, g(x)],$$
$$[n;\,a(x);\, b(x)]\cdot [m;\, f(x)]=[nm;\, f(x)],$$
$$[m;\, f(x)]\cdot [n;\, a(x);\, b(x)]=[nm;\, a(x)\cdot (x-1)^{-1}+b(x)\cdot f(x^{n})].\eqno(2.1) $$

\subsection*{3. Main theorems}

\setcounter{lemma}{0} \setcounter{theorem}{0} \setcounter{section}{3}

Suppose that $u$ is the degree of the polynomial $\psi (x)$, mentioned in
the previous section. Then $u$ is equal to the order of the element $p$ in
the group of units of the ring $Z_{q}$. The following theorem gives a
characterization of Schmidt groups by their endomorphism semigroups.

\begin{theorem}
\label{T: 3.1} Let $G$ be a finite group. Then $G$ is isomorphic to a
Schmidt group with parameters $p,\, q$ and $v$ if and only if there exists
$x\in I_{0}(G)$ such that

$1^{0}$ $K_{G}(x)\cong \hbox{End} (C(q^{v} ))$;

$2^{0}$ $H_{G}(x)=\{ 0\} $;

$3^{0}$ $I_{0}(G)=[x]$, where $[x]=\{ \, y\in \hbox{End} (G)\ \vert \
y^{2}=y,\ xy=y,\ yx=x\, \} $;

$4^{0}$ $\vert I_{0}(G)\vert =p^{u}$;

$5^{0}$ $\hbox{End} (G)\setminus \hbox{Aut} (G)=\cup _{y\in I_{0}(G)}
K_{G}(y)$;

$6^{0}$ $z\in \cap _{y\in I_{0}(G)} K_{G}(y)$ if and only if $z^{v}=0$;

$7^{0}$ $D_{G}(x)$ is a $p'$-subgroup of $V_{G}(x)$;

$8^{0}$ each Sylow $p$-subgroup of $V_{G}(x)$ is an elementary abelian
group of the

\hskip0.5cm order $p^{u}$.
\end{theorem}

{\bf Proof of the necessity of Theorem 3.1.}

Assume that $G$ is a Schmidt group with parameters $p,\, q$ and $v$. Then
$$G=G'\leftthreetimes <b>,\eqno(3.1) $$
where $<b>\cong C(q^{v})$. Denote by $x$ the projection of $G$ onto its
subgroup $<b>$. By lemma 2.1 and $\hbox{Im} \, x=<b>$, property $1^{0}$
holds.

Suppose that $z\in H_{G}(x)$. Then $zx=0$ and $xz=z$. Hence, $\hbox{Im} \,
z\subset \hbox{Ker} \, x=G'\subset \hbox{Ker} \, z$ and $bz\in G'$. Since
$b$ is a $q$-element and $G'$ is a $p$-subgroup of $G$, then $bz=1$ and
$Gz=(G'\leftthreetimes <b>)z=<1>$. Therefore, $z=0$, $H_{G}(x)=\{ 0\} $ and
property $2^{0}$ is true.

As the monomorphism $^{*}:\hbox{End} (G)\longrightarrow \hbox{End}
(G^{*})$, defined in the previous section, induces an isomorphism between
the semigroups of all proper endomorphisms of $G$ and $G^{*}=G/G''$, then
it is sufficient to prove properties $3^{0}-6^{0}$ for $G=G^{*}$ and
$x=x^{*}$.

By (2.1), $[0;\, f(x)]$ is the zero of the semigroup $\hbox{End} (G^{*})$
and
$$I_{0}(G^{*})=\{ \, [1;\, f(x)]\ \vert \ f(x)\in \overline{Z} _p[x]\, \}
.\eqno(3.2) $$ Since $[1;\, f(x)]\cdot [1;\, g(x)]=[1;\, g(x)]$, then
$I_{0}(G^{*})=[x^{*}]$ and property $3^{0}$ is true.

The idempotents $[1;\, f(x)]$ and $[1;\, g(x)]$ of $\hbox{End} (G^{*})$ are
equal if and only if $f(x)=g(x)$. Hence, $\vert I_{0}(G^{*}\vert =\vert
\overline{Z} _p[x]\vert =p^{u}$ and property $4^{0}$ holds.

It follows from (2.1) that
$$[n;\, f(x)]\cdot [1;\, f(x)]=[1;\, f(x)]\cdot [n;\, f(x)]=[n;\, f(x)],$$
i.e. $[n;\, f(x)]\in K_{G^{*}}([1;\, f(x)])$ and, therefore, property
$5^{0}$ is true.

Assume that
$$[n;\, g(x)]\in \cap _{y\in I_{0}(G^{*})}K_{G^{*}}(y).$$
By (3.2), $[n;\, g(x)]\in K_{G^{*}}([1;\, f(x)])$ for each $f(x)\in
\overline{Z} _p[x]$, i.e.
$$[n;\, g(x)]\cdot [1;\, f(x)]=[1;\, f(x)]\cdot [n;\, g(x)]=[n;\, g(x)].$$
Since the first product in the last equalities is equal to $[n;\, f(x)]$,
then $[n;\, g(x)]$ $=[n;\, f(x)]$ for each $f(x)\in \overline{Z} _p[x]$.
Hence, $[n;\, g(x)]=[n;\, 0]$ and $q$ is a divisor of $n$. Conversely, if
$q$ is a divisor of $n$, then $[n;\, 0]\in \cap _{y\in
I_{0}(G^{*})}K_{G^{*}}(y)$. Therefore,
$$\cap _{y\in I_{0}(G^{*})}K_{G^{*}}(y)=\{ \, [n;\, 0]\ \vert \ n\in q\cdot
Z_{q^{v}}\, \} .\eqno(3.3) $$
It follows from (3.3) that
$$z=[n;\, f(x)]\in \cap _{y\in I_{0}(G^{*})}K_{G^{*}}(y)$$
if and only if $z^{v}=[n^{v};\, f(x)]=0=[0;\, f(x)]$. Property $6^{0}$ is
proved.

For the proof of properties $7^{0}$ and $8^{0}$ we find first $\vert
V_{G^{*}}(x^{*})\vert $ and $\vert D_{G^{*}}(x^{*})\vert $.

Assume that $z\in \hbox{Aut} (G^{*})$. Then $z=[n;\, a(x);\, b(x)]$ for
some $n\in Z(q^{v});\ a(x),\, b(x)\in \overline{Z} _p[x]$, where $b(x)\ne
0$ and $q$ is not a divisor of $n$. It was proved in \cite{pP01} that
$x^{*}=[1;\, 0]$. The automorphism $z$ of $G^{*}$ belongs to
$V_{G^{*}}(x^{*})$ if and only if $z\cdot x^{*}=x^{*}$, i.e.,
$$[n;\, a(x);\, b(x)]\cdot [1;\, 0]=[n;\, 0]=[1;\, 0].$$
Hence, $n=1$ and
$$V_{G^{*}}(x^{*})=\{ \, [1;\, a(x);\, b(x)]\ \vert \ a(x),\, b(x)\in \overline{Z}
_p[x];\ b(x)\ne 0\, \} .$$ Since $\vert \overline{Z} _p[x]\vert =p^{u}$,
then
$$\vert V_{G^{*}}(x^{*})\vert =p^{u}(p^{u}-1).\eqno(3.4) $$

Let us calculate $\vert D_{G^{*}}(x^{*})\vert $. The group
$D_{G^{*}}(x^{*})$ consists of $z=[1;\, a(x);$ $\, b(x)]\in
V_{G^{*}}(x^{*})$ such that $x^{*}\cdot z=x^{*}$, i.e.,
$$[1;\, 0][1;\, a(x);\, b(x)]=[1;\, \frac{a(x)}{x-1} ]=[1;\, 0].$$
Hence, $a(x)=0$ and
$$\vert D_{G^{*}}(x^{*})\vert =p^{u}-1.\eqno(3.5) $$
Clearly, $(D_{G}(x))^{*}\subset D_{G^{*}}(x^{*})$, and, therefore, property
$7^{0}$ is true.

As $G''\subset Z(G),\ Z(G/G'')$ is a $q$-group and $G'/G''$ is an
elementary abelian $p$-group of the order $p^{u}$, then $G''\subset G'\cap
Z(G)$ and
$$Z(G/G'')\cap (G'/G'')=\{ 1\} .\eqno(3.6) $$
If $g\in G'\cap Z(G)$, then $gG''\in Z(G/G'')\cap (G'/G'')$ and, in view of
(3.6), $gG''=G''$ and $g\in G''$. Hence, $G'\cap Z(G)\subset G'',\
G''=G'\cap Z(G)$ and
$$\widehat{G'} =\{ \, \hat{g} \ \vert \ g\in G'\, \}
\cong G'/(G'\cap Z(G))=G'/G''.$$ Therefore, $\widehat{G'} $ and
$(\widehat{G'})^{*}$ are elementary abelian groups of the order $p^{u}$. By
lemma 2.2, $\widehat{G'} \subset V_{G}(x)$. Hence, $(\widehat{G'})^{*}
\subset V_{G^{*}}(x^{*})$ and, by (3.4), $(\widehat{G'})^{*} $ is a Sylow
$p$-subgroup of $V_{G^{*}}(x^{*})$. Then $\widehat{G'} $ is a Sylow
$p$-subgroup of $V_{G}(x)$ and property $8^{0}$ holds.

The necessity of theorem 3.1 is proved.

Assume now that $G$ is a finite group, $x\in I_{0}(G)$ and $x$ satisfies
properties $1^{0}$-$8^{0}$, formulated in theorem 3.1. Under these
assumptions we shall prove the following lemmas.

\begin{lemma}
\label{L: 3.1} Let $S$ be an arbitrary Sylow $p$-subgroup of $G$. Then
there exists $b\in G$, satisfying the following conditions:
$$G=\hbox{Ker} \, x\leftthreetimes \hbox{Im} \, x,\ \hbox{Im} \, x=<b>\cong
C(q^{v}),\eqno(3.7) $$
$$I_{0}(G)=[x]=\{ \, x\hat{g} \ \vert \ g\in G\, \} ,\eqno(3.8) $$
$$I_{0}(G)=\{ \, x\hat{g} \ \vert \ g\in \hbox{Ker} \, x\, \} ,\eqno(3.9) $$
$$[\hbox{Ker} \, x:C_{\hbox{\scriptsize Ker} \, x}(b)]=p^{u}=[S:C_{S}(b)],\eqno(3.10) $$
$$I_{0}(G)=\{ \, x\hat{s} \ \vert \ s\in S\, \} .\eqno(3.11) $$
\end{lemma}

{\bf Proof.} Since $x$ is an idempotent of $\hbox{End} (G)$, then
$G=\hbox{Ker}\, x\leftthreetimes \hbox{Im}\, x$ and $K_{G}(x)\cong
\hbox{End} (\hbox{Im} \, x)$. By $1^{0}$, $\hbox{End} (\hbox{Im} \, x)\cong
\hbox{End} (C(q^{v}))$. As each finite abelian group is determined by its
endomorphism semigroup in the class of all groups (\cite{pP75} , Theorem
4.2), then $\hbox{Im}\, x\cong C(q^{v})$ and there exists $b\in G$ such
that $\hbox{Im}\, x=<b>\cong C(q^{v})$.

By property $2^{0}$, $\hbox{Ker}\, x$ is a $q'$-group. Indeed, otherwise
there exists an element $g\in \hbox{Ker}\, x$ of the order $q$ and a
non-zero $z\in H_{G}(x)$:
$$(\hbox{Ker}\, x)z=\{ 1\} ,\ bz=g.$$
Therefore, $<b>$ is a Sylow $q$-subgroup of $G$. Since $\hbox{Ker}\, x\lhd
G$ and all Sylow subgroups related to a prime are conjugate, then all
$q'$-elements belong to $\hbox{Ker}\, x$.

By property $3^{0}$, $I_{0}(G)=[x]$. Therefore, for the proof of (3.8) and
(3.9) it is sufficient to prove the inclusions
$$[x]\subset \{ \, x\hat{g} \ \vert \ g\in \hbox{Ker}\, x \, \} \ \,
\hbox{and} \ \, \{ \, x\hat{g} \ \vert \ g\in G\, \} \subset I_{0}(G).$$

Suppose $h,\, g\in G$. As $\hbox{Im} \, x$ is abelian, then the commutator
$[hx,\, g]$ belongs to $\hbox{Ker}\, x$ and
$$h(x\hat{g} )^{2}=(g^{-1}\cdot hx\cdot g)(x\hat{g} )=(hx\cdot [hx,\,
g])(x\hat{g})=$$
$$=(hx^{2})\hat{g} =h(x\hat{g} ),$$
i.e. $(x\hat{g} )^{2}=x\hat{g} \in I_{0}(G)$. Hence, $\{ \, x\hat{g} \
\vert \ g\in G\, \} \subset I_{0}(G)$.

If $y\in [x]$, then $\hbox{Ker}\, x=\hbox{Ker}\, y$ and $\hbox{Im}\, x\cong
\hbox{Im}\, y$. As $\hbox{Im}\, x$ is a Sylow $q$-subgroup of $G$, then so
is $\hbox{Im}\, y$ and there exists $g\in \hbox{Ker}\, x$ such that
$\hbox{Im}\, y=(\hbox{Im}\, x)\hat{g} =G(x\hat{g} )=\hbox{Im}\, (x\hat{g}
)$. Since $x\hat{g} \in I_{0}(G),\ \hbox{Im}\, y=\hbox{Im}\, (x\hat{g} )$
and $\hbox{Ker}\, y=\hbox{Ker}\, x=\hbox{Ker}\, (x\hat{g} )$, then
$y=x\hat{g} $. Consequently, $[x]\subset \{ \, x\hat{g} \ \vert \ g\in
\hbox{Ker}\, x \, \}  $. Equalities (3.8) and (3.9) are proved.

In view of (3.9), $\hbox{Im}\, x=<b>$ and property $4^{0}$, the first
equality of (3.10) follows. For the proof of the second equality of (3.10)
remark that $S\subset \hbox{Ker}\, x$. Clearly,
$$[S:C_{S}(b)]\leq [\hbox{Ker}\, x:C_{\hbox{\scriptsize Ker} \, x}
(b)]=p^{u}.$$ Assume that $S_{0}$ is a Sylow $p$-subgroup of
$C_{\hbox{\scriptsize Ker} \, x} (b)$ such that $C_{S}(b)\subset S_{0}$.
Then $\vert C_{\hbox{\scriptsize Ker} \, x} (b)\vert =\vert S_{0}\vert
\cdot t$, where $p\nmid t$. By the first equality of (3.10),
$$\vert \hbox{Ker}\, x\vert =p^{u}\cdot \vert C_{\hbox{\scriptsize Ker} \, x} (b)\vert
=p^{u}\cdot \vert S_{0}\vert \cdot t$$ and, therefore, $\vert S\vert
=p^{u}\cdot \vert S_{0}\vert $. If $[S:C_{S}(b)]<p^{u}$, then $\vert
C_{S}(b)\vert >\vert S\vert :p^{u}=\vert S_{0}\vert $. This contradicts the
inclusion $C_{S}(b)\subset S_{0}$. Hence, $[S:C_{S}(b)]=p^{u}$ and
equalities (3.10) hold.

Now equality (3.11) follows already from (3.9) and (3.10). The lemma is
proved.

In the next part of this section it will be assumed that $S$ is an arbitrary
Sylow $p$-subgroup of $G$ and $\hbox{Im}\, x=<b>$.

\begin{lemma}
\label{L: 3.2} The element $b$ satisfies the following two properties
$$C_{S}(b)\subset Z(G),\eqno(3.12) $$
$$b^{q}\in Z(G).\eqno(3.13) $$
\end{lemma}

{\bf Proof.} By lemma 2.3, $\widehat{C_{S}(b)} =\{ \, \hat{s} \ \vert \ s\in
C_{S}(b)\, \} \subset D_{G}(x)$. Clearly, $\widehat{C_{S}(b)} $ is a
$p$-group. On the other hand, by property $7^{0}$, $\widehat{C_{S}(b)} $ is
a $p'$-group. Hence, $\widehat{C_{S}(b)} =\{ 1\} $ and $C_{S}(b)\subset
Z(G)$. Inclusion (3.12) is proved.

Choose $g\in G$ and denote $y=x\hat{g} $. By (3.8), $y\in I_{0}(G)$. Define
now a map $z:G\longrightarrow G$ as follows:
$$(cd)z=c^{q};\ \, c\in \hbox{Im}\, y=<g^{-1}bg>,\ d\in \hbox{Ker}\, y.$$
As $G=\hbox{Ker}\, y\leftthreetimes \hbox{Im}\, y$, then $z$ is defined
everywhere on $G$. It is easy to check that $z$ is an endomorphism of $G$
and $z^{v}=0$. By property $6^{0}$, $z\in \cap _{u\in I_{0}(G)} K_{G}(u)$.
Therefore, $z\in K_{G}(x)$ and $zx=xz$. In view of lemma 2.4, $(\hbox{Im}\,
x)z\subset \hbox{Im}\, x=<b>$ and $bz=b^{r}$ for some integer $r$. Since
the group $\hbox{Im}\, z=<g^{-1}bg>\cong G/\hbox{Ker}\, z$ is abelian, then
$G'\subset \hbox{Ker}\, z$ and $[b,\, g]z=1$. Hence,
$$b^{r}=bz=(bz)([b,\, g]z)=(b[b,\, g])z=$$
$$=(g^{-1}bg)z=g^{-1}b^{q}g=b^{q}[b^{q},\, g]$$
and
$$[b^{q},\, g]=b^{r-q}\in \hbox{Im}\, x\cap G'.$$
As $\hbox{Im}\, x\cong G/\hbox{Ker}\, x$ is abelian, then $G'\subset
\hbox{Ker}\, x$ and $[b^{q},\, g]\in \hbox{Im}\, x\cap \hbox{Ker}\, x$. By
(3.7), $\hbox{Im}\, x\cap \hbox{Ker}\, x=\{ 1\} $. Therefore, $[b^{q},\,
g]=1$. Since $g$ is an arbitrary element of $G$, then $b^{q}\in Z(G)$.
Equality (3.13) is proved, and so is lemma 3.2.

\begin{lemma}
\label{L: 3.3} There exists a Sylow $p$-subgroup $S$ of $G$ such that $b\in
N_{G}(S)$ and $<b,\, S>=S\, \leftthreetimes <b>$.
\end{lemma}

{\bf Proof.} Suppose that $S_{1}$ is an arbitrary Sylow $p$-subgroup of $G$.
Denote $N=N_{G}(S_{1})$. By (3.13), $<b^{q}>\subset N$ and, hence,
$q^{v-1}\vert \, \vert N\vert $. Assume that $q^{v}\, \nmid \, \vert N\vert
$. Then $<b^{q}>$ is a Sylow $q$-subgroup of $N$. Since $b^{q}\in Z(G)$,
then $N=<b^{q}>\times N_{1}$, where $N_{1}$ is a Hall $q'$-subgroup of $N$.
In view of (3.7), $N_{1}\subset \hbox{Ker} \, x$ and, therefore, $N\subset
<\hbox{Ker} \, x,\, b^{q}>$. Due to \cite{dR} , Theorem 1.6.18, $<\hbox{Ker}
\, x,\, b^{q}>$ coincides with its normalizator in $G$. This contradicts
(3.7). Consequently, $q^{v}\vert \, \vert N\vert $.

Since $q^{v}\vert \, \vert N\vert $ and $<b>$ is a Sylow $q$-subgroup of
the order $q^{v}$ of $G$, then $N$ consists a Sylow $q$-subgroup of $G$ and
there exists $g\in G$ such that $g^{-1}\cdot <b>\cdot g\subset
N=N_{G}(S_{1})$. Denote $S=gS_{1}g^{-1}$. Then $<b>\subset N_{G}(S)$ and
$<b,\, S>=S\, \leftthreetimes <b>$. The lemma is proved.

In the next part of this section it will be assumed that $S$ is a Sylow
$p$-subgroup of $G$ such that $b\in N_{G}(S)$ and $<b,\, S>=S\,
\leftthreetimes <b>$.

\begin{lemma}
\label{L: 3.4} The subgroup $S$ of $G$ splits up as follows
$$S=\{ \, [b,\, s]\ \vert \ s\in S\, \} \cdot C_{S}(b).\eqno(3.14) $$
\end{lemma}

{\bf Proof.} If $s\in S$ then $b^{-1}s^{-1}b\in S$ and $[b,\,
s]=b^{-1}s^{-1}b\cdot s\in S$. Assume that $[b,\, s_{1}]$ and $[b,\,
s_{2}]$ belong to a common left coset of $C_{S}(b)$ in $S$ ($s_{1},\,
s_{2}\in S$). Then $[b,\, s_{1}]=[b,\, s_{2}]\cdot c$ for some $c\in
C_{S}(b)\subset Z(G)$ and
$$s_{1}^{-1}bs_{1}=b\cdot [b,\, s_{1}]=b\cdot [b,\, s_{2}]\cdot
c=s_{2}^{-1}bs_{2}\cdot c.$$ As the order of $s_{1}^{-1}bs_{1}$ and
$s_{2}^{-1}bs_{2}$ is $q^{v}$ and $c$ is a $p$-element which belongs to
$Z(G)$, then $c=1,\ [b,\, s_{1}]=[b,\, s_{2}]$ and
$s_{1}^{-1}bs_{1}=s_{2}^{-1}bs_{2}$, i.e., $s_{1}$ and $s_{2}$ belong to a
common coset of $C_{S}(b)$ in $S$. Hence, different elements of the form
$[b,\, s]$ ($s\in S$) belong to different left cosets of $C_{S}(b)$ in $S$
and the number of elements of this form is equal to the number of cosets of
$C_{S}(b)$ in $S$. Consequently, equality (3.14) holds. The lemma is proved.

\begin{lemma}
\label{L: 3.5} The group $G$ is a semidirect product of its subgroups $S$
and $<b>$:
$$G=S\, \leftthreetimes <b>.\eqno(3.15) $$
\end{lemma}

{\bf Proof.} First we shall prove that
$$S\, \leftthreetimes <b>\lhd \, G.\eqno(3.16) $$
In view of (3.12) and (3.14), it is necessary to show that $g^{-1}bg,\,
g^{-1}[b,\, s]g\in <b,\, S>$ for each $g\in G,\ s\in S$. By (3.8) and
(3.11), $x\hat{g} =x\hat{s_{1}} $ and $x\widehat{sg} =x\widehat{s_{2}} $
for some $s_{1},\, s_{2}\in S$. On the other hand, $bx=b$. Hence,
$$g^{-1}bg=b\hat{g} =b(x\hat{g} )=b(x\hat{s_{1}} )=b\hat{s_{1}}
=s_{1}^{-1}bs_{1}\in <b,\, S>,$$
$$b\widehat{sg} =b(x\widehat{sg} )=b(x\widehat{s_{2}} )=b\widehat{s_{2}}
,$$
$$g^{-1}[b,\, s]g=g^{-1}b^{-1}s^{-1}bsg=g^{-1}b^{-1}g(b\widehat{sg} )=$$
$$=(g^{-1}bg)^{-1}\cdot (b\widehat{s_{2}})=(g^{-1}bg)^{-1}(s_{2}^{-1}bs_{2})
\in <b,\, S>.$$ Equality (3.16) is proved.

Let us prove now (3.15). By (3.16) and the theorem of Schur and Zassenhaus
(\cite{dR}, Theorem 9.1.2) there exists a $\{ p,\, q\} '$-subgroup $C$ of
$G$ such that $G=(S\, \leftthreetimes <b>)\leftthreetimes C$. Let $y$ be the
projection of $G$ onto $C$. Clearly, $y\neq 1$. If $y\neq 0$, then $y\in
I_{0}(G)$ and, by (3.8), $y=x\hat{g} $ for some $g\in G$, i.e. $\hbox{Im}
\, y=C=\hbox{Im} \, (x\hat{g} )\cong \hbox{Im} \, x\cong C(q^{v})$. This
contradicts to the fact that $C$ is a $\{ p,\, q\} '$-group. Therefore,
$y=0$ and $G=\hbox{Ker} \, y=S\, \leftthreetimes <b>$. The lemma is proved.

\begin{lemma}
\label{L: 3.6} Let $p$ and $q$ be different primes and $P$ be finite group
such that

\noindent a) $P$ is non-abelian;

\noindent b) $P=H\, \leftthreetimes <d>$, where $H$ is a $p$-group and
$<d>\cong C(q^{v})\ \, (v\ge 1)$;

\noindent c) $d^{q}\in Z(P)$;

\noindent d) $H$ is abelian;

\noindent e) if $h\in H\setminus <1>$, then $P=<h,\, d>$.

\noindent Then $P$ is a group of Miller and Moreno with parameters $p,\, q$
and $v$. If $P$ satisfies instead of properties a), d) and e) the properties

\noindent a') $P$ is non-nilpotent,

\noindent d') $C_{H}(d)\subset Z(P)$,

\noindent e') if $h\in H\setminus C_{H}(d)$, then $P=<h,\, d>$,

\noindent then $P$ is a Schmidt group with parameters $p,\, q$ and $v$.
\end{lemma}

{\bf Proof.} Assume that $P$ satisfies properties a)-e) ({\it the first
case}) or a'), b), c), d'), e') ({\it the second case}). Since different
Sylow $q$-subgroups of $P$ are conjugate and each $q$-element of $P$
belongs to some Sylow $q$-subgroup of $P$, then all $q$-elements of $P$
have a form $g^{-1}d^{i}g$ ($g\in P,\ i\in Z(q^{v})$).

Suppose that $N$ is a proper subgroup of $P$. Assume first that $N$ does
not contain $q$-elements of the order $q^{v}$. By c), all $q$-elements of
$N$ have a form $g^{-1}d^{iq}g=d^{iq}$ and, therefore, $N\subset <H,\,
d^{q}>=H\times <d^{q}>$. Hence, in the first case $N$ is abelian and in the
second case $N$ is nilpotent.

Assume now that $N$ contains a $q$-element of the order $q^{v}$, i.e.,
$g^{-1}d^{i}g\in N$, where $g\in P,\ i\in Z(q^{v})$ and $q\nmid i$. Then
$g^{-1}dg\in N$ and $d\in gNg^{-1}$.

Let us consider the first case. If $H\cap gNg^{-1}\ne <1>$, then, by e),
$gNg^{-1}=P,\ N=g^{-1}Pg=P$. It is impossible, because $N$ is a proper
subgroup of $P$. Therefore, $H\cap gNg^{-1}=<1>$. In view of $d\in
gNg^{-1}$ and b), $gNg^{-1}=<d>$, i.e. $N$ is abelian.

Let us now check the second case. It was showed that $d\in gNg^{-1}$. Since
all $p$-elements of $P$ belong to $H$, then $H\cap gNg^{-1}$ is a Sylow
$p$-subgroup of $gNg^{-1}$. Assume that $H\cap gNg^{-1}\not \subset
C_{H}(d)$ and choose $s\in (H\cap gNg^{-1})\setminus C_{H}(d)$. By e'),
$P=<s,\, d>$. On the other hand, $<s,\, d>\subset gNg^{-1}$. Hence,
$P=gNg^{-1}$ and $N=g^{-1}Pg=P$. It is impossible, because $N$ is a proper
subgroup of $P$. Therefore, $H\cap gNg^{-1}\subset C_{H}(d)$. By d'),
$gNg^{-1}=(H\cap gNg^{-1})\times <d>$ and $gNg^{-1}$ is nilpotent. Hence,
$N$ is nilpotent, too.

We have proved that in the first case all proper subgroups of $P$ are
abelian and in the second case all proper subgroups of $P$ are nilpotent.
Consequently, in the first case $P$ is a group of Miller and Moreno and in
the second case $P$ is a Schmidt group. The lemma is proved.

\begin{lemma}
\label{L: 3.7} Denote $\overline{G} =G/C_{S}(b),\ \overline{S}
=S/C_{S}(b),\ \bar{g} =g\cdot C_{S}(b)$ ($g\in G$). Then

\noindent 1) $\overline{S} $ is a non-trivial elementary abelian $p$-group;

\noindent 2) $\overline{G} =\overline{S} \, \leftthreetimes <\bar{b} >$;

\noindent 3) $G$ does not contain normal subgroups of the index $p$;

\noindent 4) $\bar{b} $ induces a non-trivial inner automorphism on
$\overline{S} $;

\noindent 5) if $\bar{s} \in \overline{S} \setminus <\bar{1} >$, then
$\overline{G} =<\bar{s} ,\, \bar{b} >$.
\end{lemma}

{\bf Proof.} In view of lemma 3.2, $C_{S}(b)\subset Z(G)$. Therefore,
$S\cap Z(G)=C_{S}(b)$ and $\widehat{S} \cong (S\cap Z(G))=\overline{S} $.
By lemma 2.2, $\widehat{S} \subset V_{G} (x)$. Hence, by $8^{0}$,
$\overline{S} $ is an elementary abelian $p$-group. Suppose $S=C_{S}(b)$.
Then equality (3.15) implies $G=S\times <b>$, i.e., $G$ is abelian.
Therefore, $x\hat{g} =x$ for each $g\in G$ and, by (3.8), $I_{0}(G)=\{ x\}
$. This contradicts property $4^{0}$. Consequently, $S\ne C_{S}(b)$ and
$\overline{S} \ne <1>$. Statement 1) of the lemma is proved.

Statement 2) of the lemma follows immediately from equality (3.15).

Let us prove now statement 3). By contradiction, assume that $M$ is a
normal subgroup of the index $p$ of $G$, i.e. $G/M\cong C(p)$. Choose a
generator $gM$ of $G/M$ and an element $h\in G$ of the order $p$. Then
there exists $\tau =\pi \rho \in \hbox{End} (G)$, where $\pi
:G\longrightarrow G/M$ is the natural homomorphism and $\rho
:G/M\longrightarrow <h>,\ (gM)\rho =h$, is the isomorphism. By construction,
$\tau $ is a non-zero proper endomorphism of $G$ and $\hbox{Im} \, \tau
\cong C(p)$. In view of property $5^{0}$, $\tau \in K_{G}(y)$ for some
$y\in I_{0}(G)$. Hence, $\tau y=y\tau =\tau $ and $\hbox{Im} \, \tau
\subset \hbox{Im} \, y$. By (3.8), $y=x\widehat{g_{1}} $ for some $g_{1}\in
G$. Therefore,
$$C(p)\cong \hbox{Im} \, \tau \subset \hbox{Im} \, y=\hbox{Im} (x\widehat{g_{1}}
)=(\hbox{Im} \, x)\widehat{g_{1}} \cong \hbox{Im} \, x\cong C(q^{v}).$$ It
is impossible. Consequently, $G$ does not contain normal subgroups of the
index $p$. Statement 3) is proved.

For the proof of statement 4) remark that $S$ is a normal subgroup of $G$
and $C_{S}(b)$ is invariant with respect to the automorphism $\hat{b} $.
Hence, it is correct to construct the automorphism $\mu $ of $\overline{S}
$ as follows
$$\bar{s} \mu =\overline{b^{-1}sb} =\bar{s} \widehat{\bar{b}} ,\ \, s\in
S.$$ Since $\bar{s} \mu ^{q} =\overline{b^{-q}sb^{q}} =\bar{s} $, then $\mu
^{q} =1$ and the order of $\mu $ is $q$ or $\mu =1$. If $\mu =1$, then, by
2), $\overline{G} =\overline{S} \times <\bar{b} >$ and there exists a
normal subgroup $M$ of $G$ such that $G/M\cong C(p)$ (we took into account
that $\overline{S} $ is a non-trivial elementary abelian $p$-group). This
contradicts to statement 3). Consequently, the order of $\mu $ is $q$, and
statement 4) is proved.

Let us prove now statement 5). The group $\overline{S} $ can be expanded
into the direct product
$$\overline{S} =\overline{S_{1}} \times \ldots \times \overline{S_{k}} ,$$
where $\overline{S_{1}} ,\ldots ,\, \overline{S_{k}} $ are minimal
non-trivial subgroups of $\overline{S} $ which are invariant with respect
to the automorphism $\mu $ which was defined above (\cite{dG} , Theorem
3.3.2). Clearly, $<\overline{S_{i}} ,\, \bar{b} >=\overline{S_{i}}
\leftthreetimes <\bar{b}
>$ for each $i\in \{ 1,\ldots ,\, k\} $. Since $\mu $ is a non-trivial
automorphism of $\overline{S} $ then there exists $i$ such that $\mu $ acts
on $\overline{S_{i}} $ non-identically. We can assume that $i=1$. Denote
$P=<\overline{S_{1}} ,\, \bar{b} >$. By construction and lemma 3.2, the
group $P$ satisfies properties a)-d) of lemma 3.6 (take there
$H=\overline{S_{1}} $ and $d=\bar{b} $). Since $\overline{S_{1}} $ is a
minimal non-trivial subgroup of $\overline{S} $ which is invariant with
respect to $\mu $, then $P$ satisfies also property e) of lemma 3.6. In
view of lemma 3.6, $\overline{S_{1}} \, \leftthreetimes <\bar{b} >$ is a
group of Miller and Moreno with parameters $p,\, q$ and $v$. By the
characterization of groups of Miller and Moreno with parameters $p,\, q$
and $v$ (\cite{lR47}), $\vert \overline{S_{1}} \vert =p^{u}$. As
$p^{u}=\vert \overline{S_{1}} \vert =\vert \overline{S_{1}} \vert :\vert
C_{S}(b)\vert =[S_{1}:C_{S}(b)]$, then, by (3.10), $S=S_{1},\, \overline{S}
=\overline{S_{1}} ,\, k=1$ and $\overline{G} =\overline{S} \,
\leftthreetimes <\bar{b} >=\overline{S_{1}} \, \leftthreetimes <\bar{b} >$.
Hence, $\overline{G} $ can be generated by $\bar{b} $ and an arbitrary
non-trivial element $\bar{s} \in \overline{S} $. Statement 5) is proved. The
lemma is proved.

\vskip0.2cm {\bf Proof of the sufficiency of theorem 3.1.}

\vskip0.2cm Let $G$ be a finite group. Assume that $x\in I_{0}(G)$ and $x$
satisfies properties $1^{0}-8^{0}$ of Theorem 3.1. By these assumptions,
lemmas 3.1-3.5 and 3.7 hold. Let us preserve the notations of these lemmas.
We will prove that $G$ is a Schmidt group with parameters $p,\, q$ and $v$.
For this aim we shall show that $G$ satisfies statements of lemma 3.6 (take
there $P=G,\, H=S,\, d=b$).

In view of lemmas 3.1, 3.2 and 3.5, statements b), c) and d') of lemma 3.6
hold. Statement a') of lemma 3.6 is also true. Indeed, otherwise $G$ is the
direct product of its Sylow subgroups and, therefore, $\overline{G}
=\overline{S} \times <\bar{b} >$. This contradicts statements 4) of lemma
3.7.

For the proof of property e') choose $s\in S\setminus C_{S}(b)$. Then
$\bar{s} =s\cdot C_{S}(b)$ is a non-trivial element of $\overline{S} $ and,
by property 5) of lemma 3.7, $\overline{G} =<\bar{s} ,\, \bar{b} >$ and
$G=<b,\, s,\, C_{S}(b)>$. Since $C_{S}(b)\subset Z(G)$, then  it follows
from here that $H=<b,\, s>$ is a normal subgroup of $G$. As $S$ is a
$p$-group, then $G/H$ is a $p$-group, too. If $H\ne G$, then $G/H$ is
non-trivial and there exists a normal subgroup $M$ of $G$ such that
$H\subset M$ and $G/M\cong C(p)$. This contradicts property 3) of lemma
3.7. Consequently, $G=H$ and $G$ satisfies property e') of lemma 3.6.

We have proved that $G$ satisfies properties a'), b), c), d') and e') of
lemma 3.6. By this lemma, $G$ is a Schmidt group with parameters $p,\, q$
and $v$. The sufficiency of theorem 3.1 is proved.

Theorem 3.1 is proved.

\begin{theorem}
\label{T: 3.2} Let $G$ be a Schmidt group with parameters $p,\, q,\, v$ and
$H$ be an arbitrary group such that the semigroups $\hbox{End} (G)$ and
$\hbox{End} (H)$ are isomorphic. Then the group $H$ is a Schmidt group with
parameters $p,\, q$ and $v$.
\end{theorem}

{\bf Proof.} Let $G$ and $H$ be groups as assumed in the theorem. Since $G$
is finite, then so is $H$ (\cite{jA62} , Theorem 2). By assumption, there
exists an isomorphism $T:\hbox{End} (G)\longrightarrow \hbox{End} (H)$. In
view of theorem 3.1, there exists $x\in I_{0}(G)$ which satisfies
properties $1^{0}-8^{0}$. By isomorphism $T$, the idempotent $xT\in
I_{0}(H)$ satisfies similar properties of this theorem. Applying theorem
3.1 for $H$, it follows that $H$ is a Schmidt group with parameters $p,\,
q$ and $v$. The theorem is proved.

Let $u$ be the order of $p$ in the group of units of residue-class ring
$Z_{q}$ modulo $q$. In \cite{lR56}, Proposition 3, it was proved that all
Schmidt groups with parameters $p,\, q$ and $v$ are isomorphic if and only
if $u$ is odd. In \cite {pP01}, Theorem 4.4, it was proved that the
endomorphism semigroup of the group of Miller and Moreno with parameters
$p,\, q$ and $v$ is isomorphic to the endomorphism semigroup of the Schmidt
group of the maximal order with same parameters. Therefore, by theorem 3.2,
the following theorem is true.

\begin{theorem}
\label{T: 3.3} The group of Miller and Moreno with parameters $p,\, q$ and
$v$ is determined by its endomorphism semigroup in the class of all groups
if and only if the order of $p$ in the group of units of the residue-class
ring $Z_{q}$ is odd.
\end{theorem}

 \vskip0.2cm \noindent {\bf Acknowledgment.} This work was
supported in part by the Estonian Science Foundation Research Grant 4291,
2000.

\end{document}